\documentclass[11pt]{article} 

\usepackage{tikz}
\usetikzlibrary{positioning}
\usepackage[utf8]{inputenc} 
\usepackage{amsmath,breqn}
\usepackage{amssymb}
\usepackage{graphicx}
\usepackage{epsfig,epstopdf}
\usepackage{pdflscape}
\usepackage{rotating,lscape,afterpage}

\usepackage{geometry} 
\usepackage{graphicx} 

\usepackage{booktabs} 
\usepackage{array} 
\usepackage{paralist} 
\usepackage{verbatim} 
\usepackage{subfig} 

\usepackage{fancyhdr} 
\usepackage{sectsty}
\allsectionsfont{\sffamily\mdseries\upshape} 

\usepackage[nottoc,notlof,notlot]{tocbibind} 
\usepackage[titles,subfigure]{tocloft} 

\title{Full street simplified three player Kuhn poker}
\author{John Billingham\\ School of Mathematical Sciences, \\ The University of Nottingham, \\ Nottingham NG7 2RD, UK}

\begin{document}
\begin{abstract}
We study a simplified version of full street, three player Kuhn poker, in which the weakest card, J, must be checked and/or folded by a player who holds it. The number of nontrivial betting frequencies that must be calculated is thereby reduced from 23 to 11, and all equilibrium solutions can be found analytically. In particular, there are three ranges of values of the pot size, $P$, for which there are three distinct, coexisting equilibrium solutions. We also study an ordinary differential equation model of repeated play of the game, which we expect to be at least qualitatively accurate when all players both adjust their betting frequencies sufficiently slowly and have sufficiently short memories. We find that none of the equilibrium solutions of the game is asymptotically stable as a solution of the ordinary differential equations. Depending on the pot size, the solution may be periodic, close to periodic or have long chaotic transients. In each case, the rates at which the players accumulate profit closely match those associated with  one of the equilibrium solutions of the game.
\end{abstract}
\maketitle

\section{Introduction}\label{sec_intro}
In \cite{SKP2017}, we studied a simplified version of one-third street, three player Kuhn poker that has three distinct equilibrium solutions. We also studied the dynamics of two models of repeated play of this game. One undesirable feature of the game studied in \cite{SKP2017} is that the restrictions on play are somewhat artificial. Another is that the play is over just one-third of a street. In this paper, we study a different simplification of three player Kuhn poker that retains more of the complexity of the full game, but is simple enough that its equilibrium solutions can be found analytically when it is played as a full street game. The simplification is that, of the four cards in the deck, only A, K and Q are active; a J must be checked and/or folded by a player who holds it\footnote{This could be enforced in practice by imposing a large penalty on any player who is revealed to have bet with J when her card is shown down.}. This means that the dynamics of the game involve sandbagging with A, bluffing with Q and bluff-catching with K\footnote{These terms will be explained below.}. This is complex enough to lead to an interesting and nontrivial equilibrium solution structure. In the full game, bluffing can be with either J or Q and bluff-catching with either Q or K, and overcalling with K may be profitable, which makes the task of finding the equilibrium solutions significantly more challenging.

We begin in Section~\ref{sec_rules} where we outline the rules of the game and describe the eleven nontrivial betting frequencies controlled by the players. In Section~\ref{sec_fullstreet} we derive the equilibrium solutions for all $P\geq2$. Finally, in Section~\ref{sec_odemodels}, we study an ordinary differential equation model of repeated play of the game, and find that the dynamic behaviour is strongly dependent on the pot size, $P$, and never leads to a steady, equilibrium solution.

\section{Simplified three player Kuhn poker}\label{sec_rules}
The deck contains four cards, A $>$ K $>$ Q $>$ J, and each of the three players is dealt a single card at random. There is a pot of $P\geq 2$ units. If $P< 2$ there is no incentive for any player to bet or call with worse than A, so we exclude this case from the outset. The first decision is made by Player 1, who must choose whether to check or bet one unit. If Player 1 bets, Players 2 and 3 must decide in turn whether to call or fold. If Player 1 checks, the decision passes to Player 2, and play proceeds in the same way. If the action is check/check/check or a bet is called, there is a showdown of the unfolded cards at which the player with the highest ranked card wins the pot and all the bets. If two players fold, the remaining player wins the pot. The full decision tree is shown in Figure~\ref{fig_decision}. 

\afterpage{\begin{landscape}
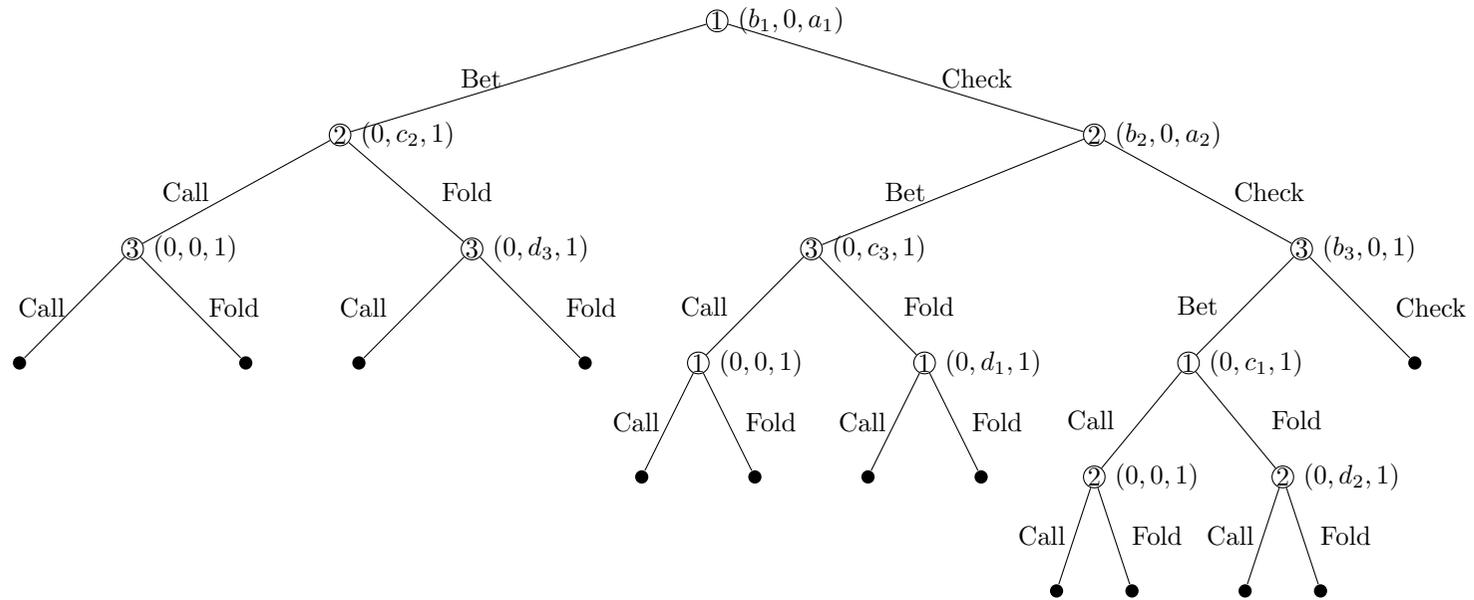
\begin{figure}
 	\begin{center}
    \small
    \begin{tikzpicture}[thin,
      level 1/.style={sibling distance=100mm},
      level 2/.style={sibling distance=55mm},
      level 3/.style={sibling distance=30mm},
      level 4/.style={sibling distance=15mm},
      level 5/.style={sibling distance=10mm},
      every circle node/.style={minimum size=1.75mm,inner sep=0mm}]

\node[circle,draw,label=right:{$(b_1,0,a_1)$}] {1}
child {
	node [circle,draw,label = right:{$(0,c_2,1)$}] {2} 
		child {
			node[circle,draw,label=right:{$(0,0,1)$}] {3}
			child {
				node [circle,fill] {} edge from parent node[left] {Call~}
				}
			child {
				node [circle,fill] {} edge from parent node[right] {~Fold}
				} edge from parent node[left] {Call~~~}
			}
		child {
			node[circle,xshift = -10mm,draw,label=right:{$(0,d_3,1)$}] {3}
			child {
				node [circle,fill] {} edge from parent node[left] {Call~~~}
				}
			child {
				node [circle,fill] {} edge from parent node[right] {~~~Fold}
				} edge from parent node[right] {~~~Fold}
			}edge from parent node[left] {Bet~~~}
	}
child {
	node [circle,draw,label=right:{$(b_2,0,a_2)$}] {2}
		child {
			node[circle,xshift = -10mm,draw,label=right:{$(0,c_3,1)$}] {3}
			child {
				node [circle,draw,label=right:{$(0,0,1)$}] {1}
			child {
				node [circle,fill] {} edge from parent node[left] {Call~}
				}
			child {
				node [circle,fill] {} edge from parent node[right] {~Fold}
				} edge from parent node[left] {Call~~~}
				}
			child {
				node [circle,draw,label=right:{$(0,d_1,1)$}] {1}
			child {
				node [circle,fill] {} edge from parent node[left] {Call~}
				}
			child {
				node [circle,fill] {} edge from parent node[right] {~Fold}
				} edge from parent node[right] {~~~Fold}
				} edge from parent node[left] {Bet~~~}
			}
		child {
			node[circle,draw,label=right:{$(b_3,0,1)$}] {3}
			child {
				node [circle,draw,label=right:{$(0,c_1,1)$}] {1}
		child {
			node[circle,draw,xshift=-5mm,label=right:{$(0,0,1)$}] {2}
			child {
				node [circle,xshift=0mm,fill] {} edge from parent node[left] {Call~}
				}
			child {
				node [circle,xshift=0mm,fill] {} edge from parent node[right] {~Fold}
				} edge from parent node[left] {Call~~~}
			}
		child {
			node[circle,draw,xshift=5mm,label=right:{$(0,d_2,1)$}] {2}
			child {
				node [circle,xshift=0mm,fill] {} edge from parent node[left] {Call~}
				}
			child {
				node [circle,xshift=0mm,fill] {} edge from parent node[right] {~Fold}
				} edge from parent node[right] {~~~Fold}
			} edge from parent node[left] {Bet~~~}
				}
			child {
				node [circle,fill] {} edge from parent node[right] {~~~Check}
				} edge from parent node[right] {~~~Check}
			} edge from parent node[right] {~~~Check}
	}
;

    \end{tikzpicture}
    \end{center}
    \caption{The decision tree for three player, full street simplified Kuhn poker. Open circles are decision nodes (labelled by the player making the decision), whilst solid circles are terminal nodes. Betting and calling frequencies with (Q, K, A) are shown at each node. A player who holds J must check and/or fold.\label{fig_decision}}
  \end{figure}
\end{landscape}}

The crucial difference between this game and standard Kuhn poker is that J is a dead card. A player who is dealt J must check and/or fold at each of their decision points. In this way, the game is reduced to a form more similar to two player Kuhn poker. The three players must decide whether to: (i) bet with A, hoping to be called by K, or check and hope a player with Q bluffs so that they can call and win an extra bet (known as sandbagging); (ii) bet with Q and hope both that A has not been dealt and that the player who holds K folds (bluffing); (iii) call a bet with K hoping that the bet was a bluff with Q (bluff-catching). The card J is in the deck simply to make the game nontrivial, since a three card deck always deals a winning A to one of the players. There are just eleven nontrivial betting frequencies that must be chosen by the players. As we shall see, we are able to determine all of the equilibrium solutions analytically. In contrast, there are 23 nontrivial betting frequencies in standard three player Kuhn poker, in which bluffing can be with J or Q, bluffcatching with Q or K and a player may wish to overcall with K, which makes the equilibria of the game much harder to find analytically. The only published equilibrium solution is for $P=3$, described in \cite{Szafron:2013:PFE:2484920.2484962}.

Note that in our simplified game: (i) betting with K is dominated by checking for every player, because two thirds of the time one of the remaining players will hold A and call; (ii) after a bet and a call, the remaining player will call only with A. Our notation for the non-trivial betting frequencies of Player $i$ is $b_i$ for betting with Q (her bluffing frequency), $c_i$ for calling with K after a bet, $d_i$ for calling with K after a bet and a fold, and $a_i$ for betting with A\footnote{$1-a_i$ is her sandbagging frequency.}. After Players 1 and 2 check, betting with A dominates checking for Player 3, so $a_3=1$. There are therefore eleven frequencies to be determined.
\newpage
\section{Equilibrium solutions}\label{sec_fullstreet}
By considering each of the 24 distinct deals of three of the four cards and the probabilities of each outcome, we find that the expected profits of the three players, $E_i$ for $i = 1$, $2$ and $3$, are given by
\begin{multline}
24 E_1 = \left\{-2 b_2 + 2 c_2 - 2 b_3 + 2 d_3 - b_2 c_3\right\}a_1 + \left\{2P-4 -\left(P+1\right)\left(c_2+d_3\right)\right\}b_1 \\ +\left\{b_2-2 + \left(P+a_2\right) b_3\right)\}c_1 + \left\{\left(P+1\right)b_2 - 2 a_2\right\} d_1 + \left(2-P\right) b_3 + \left(2+c_3 -P\right)b_2,
\label{eq_E1}
\end{multline}
\begin{multline}
24 E_2 = \left\{2 d_1 + 2 c_3 - 2 b_3 - c_1 b_3 - b_1 c_3\right\} a_2 + \left\{2P-4+2 a_1 -\left(P+1\right)\left(c_3+d_1\right)\right\} b_2\\+ \left\{ P b_1-2 a_1 \right\} c_2 + \left\{ \left(P+1\right) b_3 -2 +b_1 \right\} d_2 + \left(2-P + c_1\right) b_3 + \left(2-P\right) b_1,
\label{eq_E2}
\end{multline}
\begin{multline}
24 E_3 = \left\{2P-4+2 a_1+2a_2 - \left(P+1\right) \left(c_1 +d_2\right) \right\} b_3 + \left\{\left(P+a_1\right) b_2 -\left(2-b_1\right)a_2 \right\} c_3 + \\ \left\{ \left(P+1\right) b_1-2 a_1 \right\} d_3+ \left(2-P-c_1 \right) b_2 + 2 c_1 + \left(2-P + c_2 -d_2 \right) b_1 + 2 d_2.
\label{eq_E3}
\end{multline}
Since Player $i$ chooses $a_i$, $b_i$, $c_i$ and $d_i$ to maximise $E_i$, the eleven conditions that must be satisfied by an equilibrium solution are:
\begin{enumerate}
\item
\begin{enumerate}
\item $b_1 > 2 a_1/P$ \& $c_2 = 1$,
\item $b_1 < 2 a_1/P$ \& $c_2 = 0$,
\item $b_1 = 2 a_1/P$.
\end{enumerate}
\item
\begin{enumerate}
\item $b_1 > 2 a_1/(P+1)$ \& $d_3= 1$,
\item $b_1 < 2 a_1/(P+1)$ \& $d_3 = 0$,
\item $b_1 = 2 a_1/(P+1)$.
\end{enumerate}
\item
\begin{enumerate}
\item $c_2+d_3 < (2P-4)/(P+1)$ \& $b_1= 1$,
\item $c_2+d_3 > (2P-4)/(P+1)$ \& $b_1 = 0$,
\item $c_2+d_3 = (2P-4)/(P+1)$.
\end{enumerate}
\item
\begin{enumerate}
\item $c_2+d_3 > b_3+(1+\frac{1}{2}c_3)b_2$ \& $a_1= 1$,
\item $c_2+d_3 < b_3+(1+\frac{1}{2}c_3)b_2$ \& $a_1 = 0$,
\item $c_2+d_3 = b_3+(1+\frac{1}{2}c_3)b_2$ .
\end{enumerate}
\item
\begin{enumerate}
\item $b_2 > (2-b_1)a_2/(P+a_1)$ \& $c_3= 1$,
\item $b_2 < (2-b_1)a_2/(P+a_1)$ \& $c_3 = 0$,
\item $b_2 = (2-b_1)a_2/(P+a_1)$.
\end{enumerate}
\item
\begin{enumerate}
\item $b_2 > 2 a_2/(P+1)$ \& $d_1= 1$,
\item $b_2 < 2 a_2/(P+1)$ \& $d_1 = 0$,
\item $b_2 = 2 a_2/(P+1)$.
\end{enumerate}
\item
\begin{enumerate}
\item $c_3+d_1 < (2P-4+2a_1)/(P+1)$ \& $b_2= 1$,
\item $c_3+d_1 > (2P-4+2a_1)/(P+1)$ \& $b_2 = 0$,
\item $c_3+d_1 = (2P-4+2a_1)/(P+1)$.
\end{enumerate}
\item
\begin{enumerate}
\item $c_3+d_1 >\frac{1}{2}c_3 b_1+(1+\frac{1}{2}c_1)b_3$ \& $a_2= 1$,
\item $c_3+d_1 <\frac{1}{2}c_3 b_1+(1+\frac{1}{2}c_1)b_3$ \& $a_2 = 0$,
\item $c_3+d_1 = \frac{1}{2}c_3 b_1+(1+\frac{1}{2}c_1)b_3$ .
\end{enumerate}
\item
\begin{enumerate}
\item $b_3 > (2-b_2)/(P+a_2)$ \& $c_1= 1$,
\item $b_3 < (2-b_2)/(P+a_2)$ \& $c_1 = 0$,
\item $b_3 = (2-b_2)/(P+a_2)$.
\end{enumerate}
\item
\begin{enumerate}
\item $b_3 > (2-b_1)/(P+1)$ \& $d_2= 1$,
\item $b_3 < (2-b_1)/(P+1)$ \& $d_2 = 0$,
\item $b_3 = (2-b_1)/(P+1)$.
\end{enumerate}
\item
\begin{enumerate}
\item $c_1+d_2 < (2P-4+2a_1+2a_2)/(P+1)$ \& $b_3= 1$,
\item $c_1+d_2 > (2P-4+2a_1+2a_2)/(P+1)$ \& $b_3 = 0$,
\item $c_1+d_2 = (2P-4+2a_1+2a_2)/(P+1)$.
\end{enumerate}
\end{enumerate}
At least one of the subconditions labelled (a), (b) and (c) must be satisfied for each of the eleven conditions. There are however several possibilities that we can eliminate.
\begin{itemize}
\item 3.(a) cannot hold. ($b_1 \not = 1$)

3.(a) $\implies b_1=1$. Then 1. $\implies c_2=1$ and 2. $\implies d_3=1$. This means that $c_2+d_3 = 2 > (2P-4)/(P+1)$, in contradiction with 3.(a).
\item 7.(a) cannot hold.  ($b_2 \not = 1$)

7.(a) $\implies$ $b_2=1$. Then 5. $\implies c_3=1$ and 6. $\implies d_1=1$. This means that $c_3+d_1 = 2 > (2P-4 + 2a_1)/(P+1)$, in contradiction with 7.(a).
\item Neither 11.(a) nor 11.(b) can hold.  ($0 < b_3 < 1$)

11.(a) $\implies$ $b_3=1$. Then 9. $\implies c_1=1$ and 10. $\implies d_2=1$. This means that $c_1+d_2 = 2 > (2P-4 + 2a_1 + 2a_2)/(P+1)$, in contradiction with 11.(a).

11.(b) $\implies b_3=0$. Then 9. $\implies c_1=0$ and 10. $\implies d_2 = 0$. This means that $c_1+d_2 = 0 < (2P-4+2a_1+2a_2)/(P+1)$, in contradiction with 11.(b).

\item Either $a_1 = b_1 = 0$ or both $a_1>0$ and $b_1>0$. 

If $a_1=0$ and $b_1>0$, then 1. and 2. show that $c_2 = d_3 = 1$, and hence $c_2+d_3 = 2 > (2P-4)/(P+1)$, in contradiction with 3. If $a_1>0$ and $b_1=0$, then 1. and 2. show that $c_2 = d_3 = 0$, and hence $c_2+d_3 = 0 < (2P-4)/(P+1)$, in contradiction with 3.

\item Either $a_2 = b_2 = 0$ or both $a_2>0$ and $b_2>0$. 

If $a_2=0$ and $b_2>0$, then 5. and 6. show that $c_3 = d_1 = 1$, and hence $c_3+d_1 = 2 > (2P-4+2a_1)/(P+1)$, in contradiction with 7. If $a_2>0$ and $b_2=0$, then 5. and 5. show that $c_3 = d_1 = 0$, and hence $c_2+d_3 = 0 < (2P-4+2a_1)/(P+1)$, in contradiction with 3.
\end{itemize}
In poker terms, we have shown that no player has a range strong enough that they should always bluff with Q, and also that any value bet with A must be balanced by bluffs with Q and vice versa.
\begin{table}
\begin{enumerate}
\item
\begin{enumerate}
\item $a_1=0$ \& $b_1=0$ \& $c_2+d_3 \geq (2P-4)/(P+1)$,
\item $b_1 = 2 a_1/P$ \& $c_2 = (P-5)/(P+1)$ \& $d_3=1$ \& $0<a_1<1$ \& $0<b_1<1$ \& $P\geq 5$,
\item$b_1 = 2 a_1/(P+1)$ \& $c_2 = 0$ \& $d_3=(2P-4)/(P+1)$ \& $0<a_1<1$ \& $0<b_1<1$ \& $P\leq 5$,
\end{enumerate}
\item
\begin{enumerate}
\item $c_2+d_3 > b_3+(1+\frac{1}{2}c_3)b_2$ \& $a_1= 1$,
\item $c_2+d_3 < b_3+(1+\frac{1}{2}c_3)b_2$ \& $a_1 = 0$,
\item $c_2+d_3 = b_3+(1+\frac{1}{2}c_3)b_2$ .
\end{enumerate}
\item
\begin{enumerate}
\item $b_2 > (2-b_1)a_2/(P+a_1)$ \& $c_3= 1$,
\item $b_2 < (2-b_1)a_2/(P+a_1)$ \& $c_3 = 0$,
\item $b_2 = (2-b_1)a_2/(P+a_1)$.
\end{enumerate}
\item
\begin{enumerate}
\item $b_2 > 2 a_2/(P+1)$ \& $d_1= 1$,
\item $b_2 < 2 a_2/(P+1)$ \& $d_1 = 0$,
\item $b_2 = 2 a_2/(P+1)$.
\end{enumerate}
\item
\begin{enumerate}
\item $c_3+d_1 > (2P-4+2a_1)/(P+1)$ \& $b_2 = 0$,
\item $c_3+d_1 = (2P-4+2a_1)/(P+1)$.
\end{enumerate}
\item
\begin{enumerate}
\item $c_3+d_1 >\frac{1}{2}c_3 b_1+(1+\frac{1}{2}c_1)b_3$ \& $a_2= 1$,
\item $c_3+d_1 <\frac{1}{2}c_3 b_1+(1+\frac{1}{2}c_1)b_3$ \& $a_2 = 0$,
\item $c_3+d_1 = \frac{1}{2}c_3 b_1+(1+\frac{1}{2}c_1)b_3$ .
\end{enumerate}
\item
\begin{enumerate}
\item $b_3 > (2-b_2)/(P+a_2)$ \& $c_1= 1$,
\item $b_3 < (2-b_2)/(P+a_2)$ \& $c_1 = 0$,
\item $b_3 = (2-b_2)/(P+a_2)$.
\end{enumerate}
\item
\begin{enumerate}
\item $b_3 > (2-b_1)/(P+1)$ \& $d_2= 1$,
\item $b_3 < (2-b_1)/(P+1)$ \& $d_2 = 0$,
\item $b_3 = (2-b_1)/(P+1)$.
\end{enumerate}
\item $c_1+d_2 = (2P-4+2a_1+2a_2)/(P+1)$.
\end{enumerate}
\caption{The nine conditions that must be satisfied by an equilibrium solution. For each of conditions 1 to 8, at least one of the subconditions (a), (b) and (c) must hold.}\label{table2}
\end{table}

Finally, when $a_1>0$ and $b_1>0$, by considering 1. and 2., noting that $2a_1/P > 2 a_1/(P+1)$, we can deduce that
\[
b_1 = \frac{2a_1}{P+1},~~c_2=0,~~d_3 = \frac{2P-4}{P+1}~~\mbox{for $P\leq 5$,}\]
\begin{equation}
b_1 = \frac{2a_1}{P},~~c_2 = \frac{P-5}{P+1},~~d_3=1~~\mbox{for $P\geq 5$.}
\end{equation}
Note that a similar argument does not work immediately for $a_2>0$ and $b_2>0$ since the size of  $(2-b_1)/(P+a_1)$  relative to $2/(P+1)$ depends on $a_1$ and $b_1$. 

\begin{figure}
 \begin{center}
 \includegraphics[width=\textwidth]{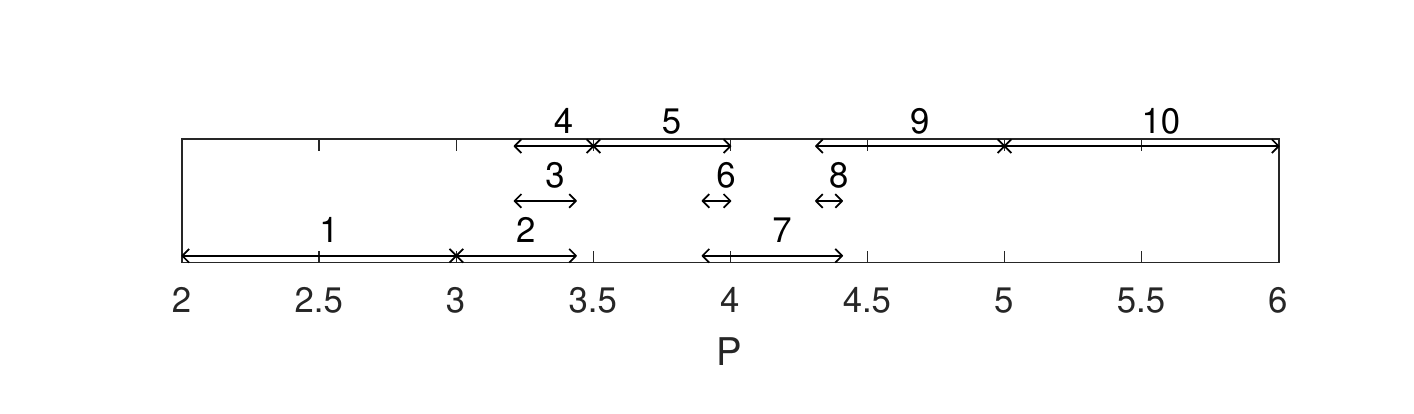}
 \caption{The ranges of validity of the ten equilibrium solutions.}\label{fig_ranges}
 \end{center}
 \end{figure}

\begin{figure}
 \begin{center}
 \includegraphics[width=0.9\textwidth]{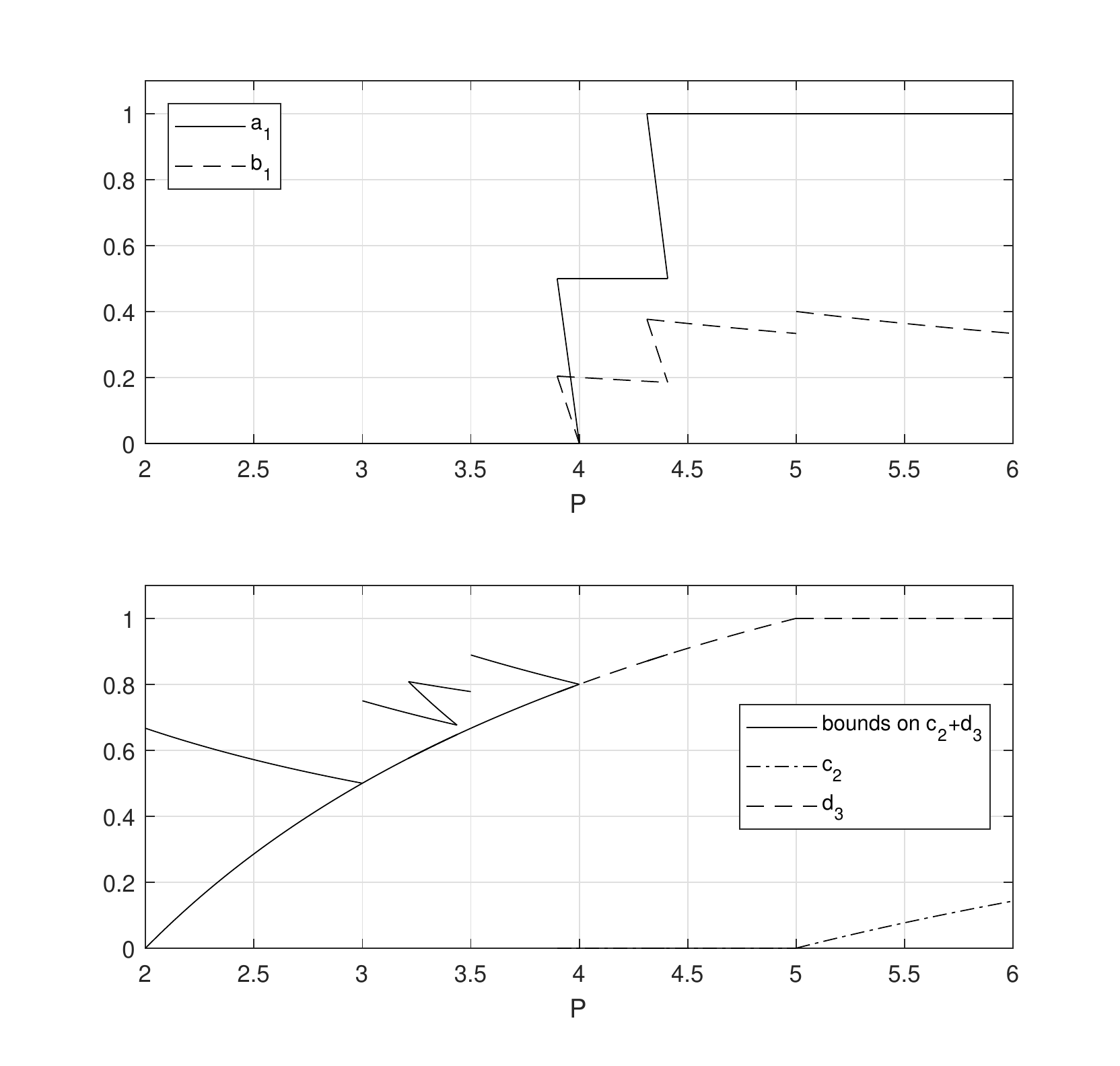}
 \caption{Player 1's equilibrium betting frequencies and the consequent calling frequencies.}\label{fig_a1b1}
 \end{center}
 \end{figure}

It is now helpful to rewrite the conditions in the more convenient form shown in Table~\ref{table2}. We used the computer algebra package Mathematica to simplify each of the $2 . 3^7 = 4374$ combinations of subcases of the nine conditions shown in Table~\ref{table2} and thereby determined all possible equilibrium solutions analytically\footnote{The code and output has been deposited as ancillary material on arXiv.}. The calculation took about 25 minutes on a single core of a standard desktop computer and the solutions found were then checked by hand. In total we find that there are ten distinct solutions that are valid in the ranges shown in Figure~\ref{fig_ranges}, along with specific solutions when $P = 2$, $3$, $7/2$ and $5$, which we label separately. 
This illustrates the fact that there are three ranges of values of the pot size, $P$, for which three equilibrium solutions coexist. 
\newpage
The equilibrium solutions are:
\[ \mbox{\bf Solution 1a:}~~a_1=b_1=a_2=b_2=0,~~0 \leq b_3\leq \frac{2}{3},~~c_1=d_2=0,\]
\begin{equation}
c_3+d_1 \leq b_3,~~c_2+d_3 \leq b_3,~~{\bf P=2}.~~\label{eq_sol1a}
\end{equation}
\[ \mbox{\bf Solution 1:}~~a_1=b_1=a_2=b_2=0,~~b_3=\frac{2}{P+1},~~c_1=0,~~d_2=\frac{2P-4}{P+1},\]
\begin{equation}
\frac{2P-4}{P+1} \leq c_3+d_1 \leq \frac{2}{P+1},~~\frac{2P-4}{P+1}\leq c_2+d_3 \leq \frac{2}{P+1},~~{\bf 2 \leq P \leq 3}.~~\label{eq_sol1}
\end{equation}
\[ \mbox{\bf Solution 2a:}~~a_1=b_1=0,~~b_2 = \frac{1}{2}a_2,~~a_2 \leq \frac{1}{2},~~b_3=\frac{1}{2},\]
\begin{equation}
c_1=0,~~d_2=\frac{1}{2}(1+a_2),~~c_3=0,~~d_1=\frac{1}{2},~~\frac{1}{2}\leq c_2+d_3 \leq \frac{1}{2}+b_2,~~{\bf P=3}.~~\label{eq_sol2a}
\end{equation}
\[\mbox{\bf Solution 2:}~~a_1 = b_1=0,~~a_2 = \frac{1}{2},~~b_2 = \frac{1}{P+1},~~b_3 = \frac{2}{P+1},\]
\begin{equation}
c_1 = 2P-6,~~d_2 = \frac{3+6P-2P^2}{P+1},~~c_3 = 0,~~d_1 = \frac{2P-4}{P+1},\label{eq_sol2}
\end{equation}
\[~~\frac{2P-4}{P+1}\leq c_2+d_3 \leq \frac{3}{P+1},~~{\bf 3 \leq P \leq P_4\equiv \frac{1}{2}\left(3+\sqrt{15}\right) \approx 3.43}.\]
\[\mbox{\bf Solution 3:}~~a_1=b_1=0,~~a_2 = \frac{1}{2}\left(4-P^2+\sqrt{P^4-4P^3+4P^2+12P+12}\right),~~b_2 = \frac{2a_2}{P+1},\]
\begin{equation}
b_3 = \frac{2-b_2}{P+a_2},~~c_1 = \frac{2P-4+2a_2}{P+1},~~d_2=0,~~c_3=0,~~d_1 = \frac{2P-4}{P+1},\label{eq_sol3}
\end{equation}
\[\frac{2P-4}{P+1} \leq c_2+d_3 \leq b_2+b_3,~~{\bf P_3 \equiv \frac{1}{4}\left(3+\sqrt{97}\right) \approx 3.21 \leq P \leq P_4 \equiv \frac{1}{2}\left(3+\sqrt{15}\right) \approx 3.43}.\]
\[\mbox{\bf Solution 4:}~~a_1=b_1=0,~~a_2 = \frac{1}{2}(5-P),~~b_2 = \frac{5-P}{P+1},~~b_3 = \frac{4(P-2)}{3(P+1)},\]
\begin{equation}
c_1 = 1,~~d_2=0,~~c_3=0,~~d_1 = \frac{2P-4}{P+1},\label{eq_sol4}
\end{equation}
\[\frac{2P-4}{P+1}\leq c_2+d_3\leq \frac{P+7}{3(P+1)},~~{\bf P_3 \equiv \frac{1}{4}\left(3+\sqrt{97}\right) \approx 3.21 \leq P \leq \frac{7}{2}}.\]
\[\mbox{\bf Solution 5a:}~~a_1=b_1=0,~~b_2 = \frac{4}{9}a_2,~~a_2 \geq \frac{3}{4},~~b_3 = \frac{4}{9},\]
\begin{equation}
c_1 = 1,~~d_2=\frac{4}{9}a_2-\frac{1}{3},~~c_3=0,~~d_1 = \frac{2}{3},\label{eq_sol5a}
\end{equation}
\[\frac{2}{3}\leq c_2+d_3\leq \frac{4}{9}(1+a_2),~~{\bf P =\frac{7}{2}}.\]
\newpage
\[\mbox{\bf Solution 5:}~~a_1=b_1=0,~~a_2 = 1,~~b_2 = \frac{2}{P+1},~~b_3 = \frac{2}{P+1},\]
\begin{equation}
c_1 = 1,~~d_2=\frac{P-3}{P+1},~~c_3=0,~~d_1 = \frac{2P-4}{P+1},\label{eq_sol5}
\end{equation}
\[\frac{2P-4}{P+1}\leq c_2+d_3\leq \frac{4}{P+1},~~{\bf \frac{7}{2} \leq P \leq 4}.\]
\[\mbox{\bf Solution 6:}~~a_1=(4-P)(P+1),~~b_1=2(4-P),~~a_2 = 1,~~b_2 = \frac{2}{P+1},\]
\begin{equation}
b_3 = \frac{2(P-3)}{P+1},~~c_1 = 1,~~d_2=\frac{5+7P-2P^2}{P+1},~~c_3=0,~~d_1 = \frac{4+8P-2P^2}{P+1},\label{eq_sol6}
\end{equation}
\[c_2=0,~~d_3=\frac{2P-4}{P+1},~~{\bf P_6 \equiv \frac{1}{2}\left(3+\sqrt{23}\right) \approx 3.90 \leq P \leq 4}.\]
\[\mbox{\bf Solution 7:}~~a_1=\frac{1}{2},~~b_1=\frac{1}{P+1},~~a_2 = 1,~~b_2 = \frac{2}{P+1},~~b_3 = \frac{2P+1}{(P+1)^2},\]
\begin{equation}
c_1 = 1,~~d_2=\frac{P-2}{P+1},~~c_3=\frac{2P^2-6P-7}{P+1},~~d_1 = \frac{4+8P-2P^2}{P+1},~~c_2=0,\label{eq_sol7}
\end{equation}
\[d_3=\frac{2P-4}{P+1},~~{\bf P_6 \equiv \frac{1}{2}\left(3+\sqrt{23}\right) \approx 3.90 \leq P \leq P_9 \equiv \frac{1}{4}\left(7+\sqrt{113}\right)\approx 4.41}.\]
\[\mbox{\bf Solution 8:}~~a_1=\frac{1}{2}(9-2P)(P+1),~~b_1=9-2P,~~a_2 = 1,~~b_2 = \frac{2}{P+1},~~b_3 = \frac{2P-7}{P+1},\]
\begin{equation}
c_1 = 1,~~d_2=\frac{6+8P-2P^2}{P+1},~~c_3=1,~~d_1 = \frac{4+8P-2P^2}{P+1},~~c_2=0,\label{eq_sol8}
\end{equation}
\[d_3=\frac{2P-4}{P+1},~~{\bf P_8 \equiv \frac{1}{4}\left(7+\sqrt{105}\right) \approx 4.31 \leq P \leq P_9 \equiv \frac{1}{4}\left(7+\sqrt{113}\right)\approx 4.41}.\]
\[\mbox{\bf Solution 9:}~~a_1=1,~~b_1=\frac{2}{P+1},~~a_2 = 1,~~b_2 = \frac{2}{P+1},~~b_3 = \frac{2P}{(P+1)^2},\]
\begin{equation}
c_1 +d_2=\frac{2P}{P+1},~~c_3=1,~~d_1 = \frac{P-3}{P+1},~~c_2=0,\label{eq_sol9}
\end{equation}
\[d_3=\frac{2P-4}{P+1},~~{\bf P_8 \equiv \frac{1}{4}\left(7+\sqrt{105}\right) \approx 4.31 \leq P \leq 5}.\]
\[\mbox{\bf Solution 10a:}~~a_1=1,~~\frac{1}{3} \leq b_1 \leq \frac{2}{5},~~a_2 = 1,~~b_2 = \frac{1}{3},~~b_3 = \frac{5}{18},\]
\begin{equation}
c_1 +d_2=\frac{5}{3}~\mbox{when $b_1 = \dfrac{1}{3}$,}~~~~c_1=\frac{2}{3},~d_2=1~\mbox{when $b_1 > \dfrac{1}{3}$,}\label{eq_sol10a}
\end{equation}
\[c_3=1,~~d_1 = \frac{1}{3},~~c_2=0,~~d_3=1,~~{\bf P= 5}.\]
\newpage
\[\mbox{\bf Solution 10:}~~a_1=1,~~b_1=\frac{2}{P},~~a_2 = 1,~~b_2 = \frac{2}{P+1},~~b_3 = \frac{2P}{(P+1)^2},\]
\begin{equation}
c_1=\frac{P-1}{P+1},~~d_2=1,~~c_3=1,~~d_1 = \frac{P-3}{P+1},~~c_2=\frac{P-5}{P+1},\label{eq_sol10}
\end{equation}
\[d_3=1,~~{\bf P \geq 5}.\]
Figure~\ref{fig_a1b1} shows the equilibrium frequencies at which Player 1 bets and Players 2 and 3 call. For small enough P ($P<P_6 \approx 3.9$), Player 1 neither bets with A nor bluffs with Q at equilibrium. For these values of $P$, the solid lines in the lower graph show the bounds between which $c_2+d_3$ must lie at equilibrium. For values above the upper bound, Player 1 can exploit her opponents by always betting with A; for values below the lower bound, by always bluffing with Q. For larger values of $P$, Player 1 has nonzero betting and bluffing frequencies at equilibrium, including two ranges of $P$ with multiple solutions. Only Player 3 catches Player 1's bluffs for $P<5$, whilst for $P > 5$, Player 1 always calls with K and Player 2 also catches bluffs with K at a nonzero frequency ($d_3=1$, $c_2>0$).

\begin{figure}
 \begin{center}
 \includegraphics[width=0.9\textwidth]{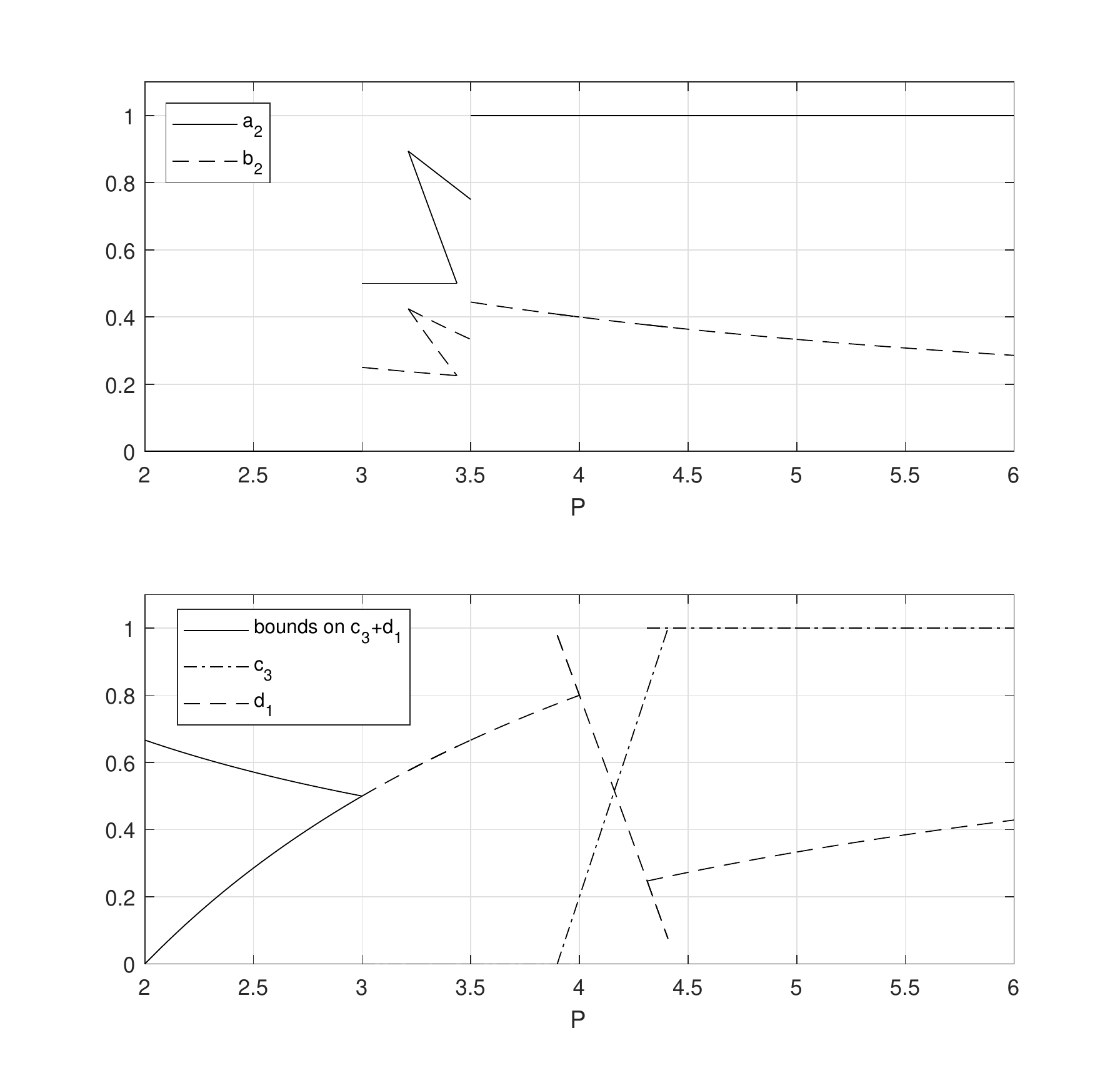}
 \caption{Player 2's equilibrium betting frequencies and the consequent calling frequencies.}\label{fig_a2b2}
 \end{center}
 \end{figure}

Figure~\ref{fig_a2b2} shows the equilibrium frequencies at which Player 2 bets and Players 3 and 1 call. For $P<3$, Player 2 does not bet at equilibrium, and there are upper and lower bounds on $c_3+d_1$.  For $P>3$, Player 2 bets at equilibrium, with one range of multiple solutions. Players 3 and 1 have calling frequencies that have a rather unintuitive structure at equilibrium. For small enough $P$, Player 1 is the bluffcatcher, but as $P$ increases, Player 3 develops a nonzero equilibrium calling frequency. For large enough $P$, Player 3 always calls with K and Player 1 has a nonzero calling frequency. This is in contrast to the frequencies shown in Figure~\ref{fig_a1b1}, where the player who acts immediately after the original bettor (Player 2) does not call with K for small enough values of $P$, but the next player to act (Player 3) does.

\begin{figure}
 \begin{center}
 \includegraphics[width=0.9\textwidth]{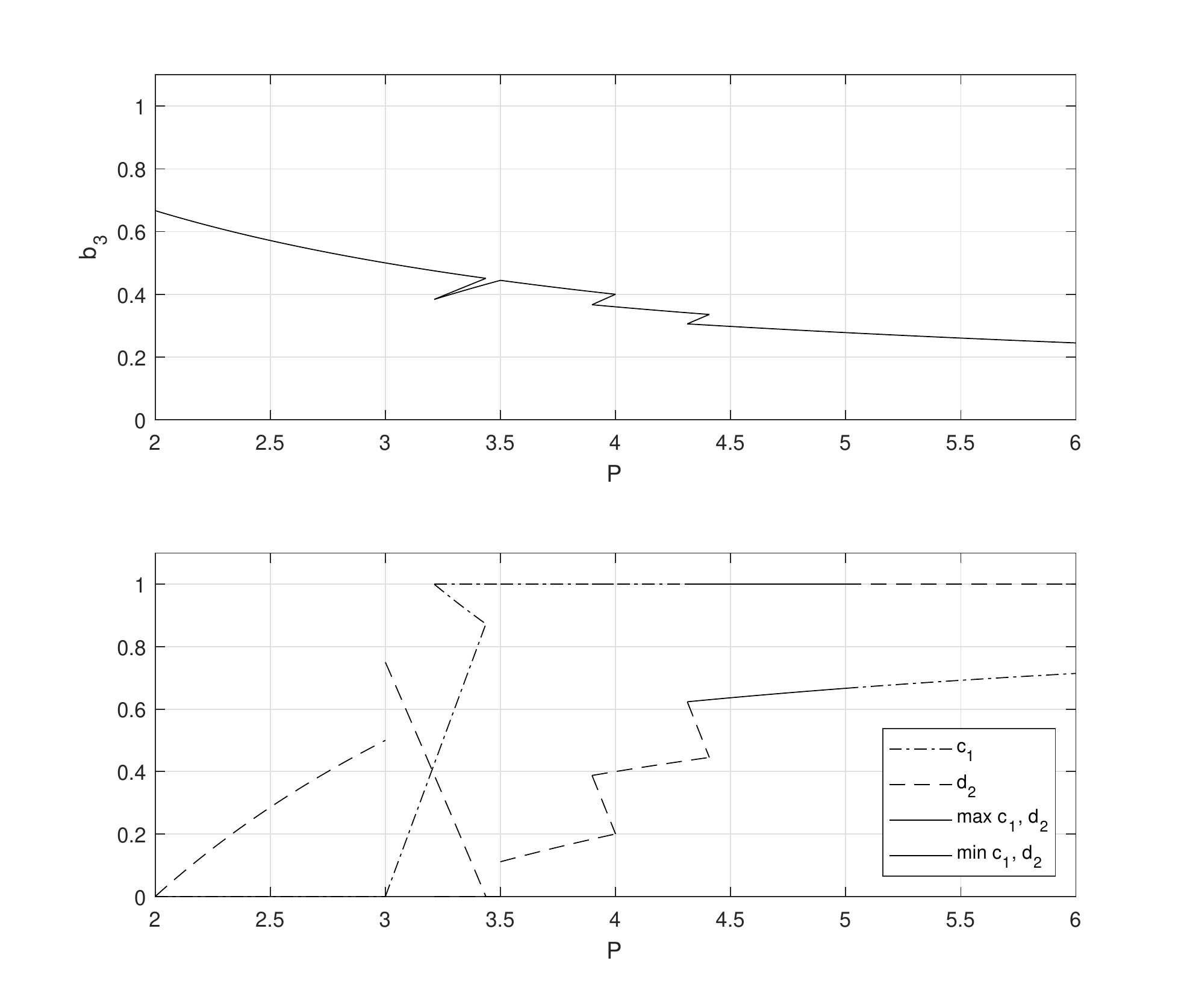}
 \caption{Player 3's equilibrium bluffing frequency and the consequent calling frequencies.}\label{fig_b3}
 \end{center}
 \end{figure}

Figure~\ref{fig_b3} shows the equilibrium frequencies at which Player 3 bluffs and Players 1 and 2 call. Since Player 3 always bets with A, she has a nonzero bluffing frequency for all $P$, which generally decreases as $P$ increases, with some irregular behaviour in the regions where three equilibrium solutions exist. The most notable feature of the calling frequencies of Players 1 and 2 is that there is a set of values of $P$ for which Player 1 always calls with K ($c_1=1$), followed by a region where there is a range of possible values of $c_1$ and $d_1$, bounded by the solid lines, followed for large enough $P$ by a set of values for which Player 2 always calls with K ($d_2=1$). It is hard to interpret these equilibrium solutions, in particular the peculiar way in which the calling frequencies vary with $P$, either mathematically or in terms of poker concepts.

\begin{figure}
 \begin{center}
 \includegraphics[width=\textwidth]{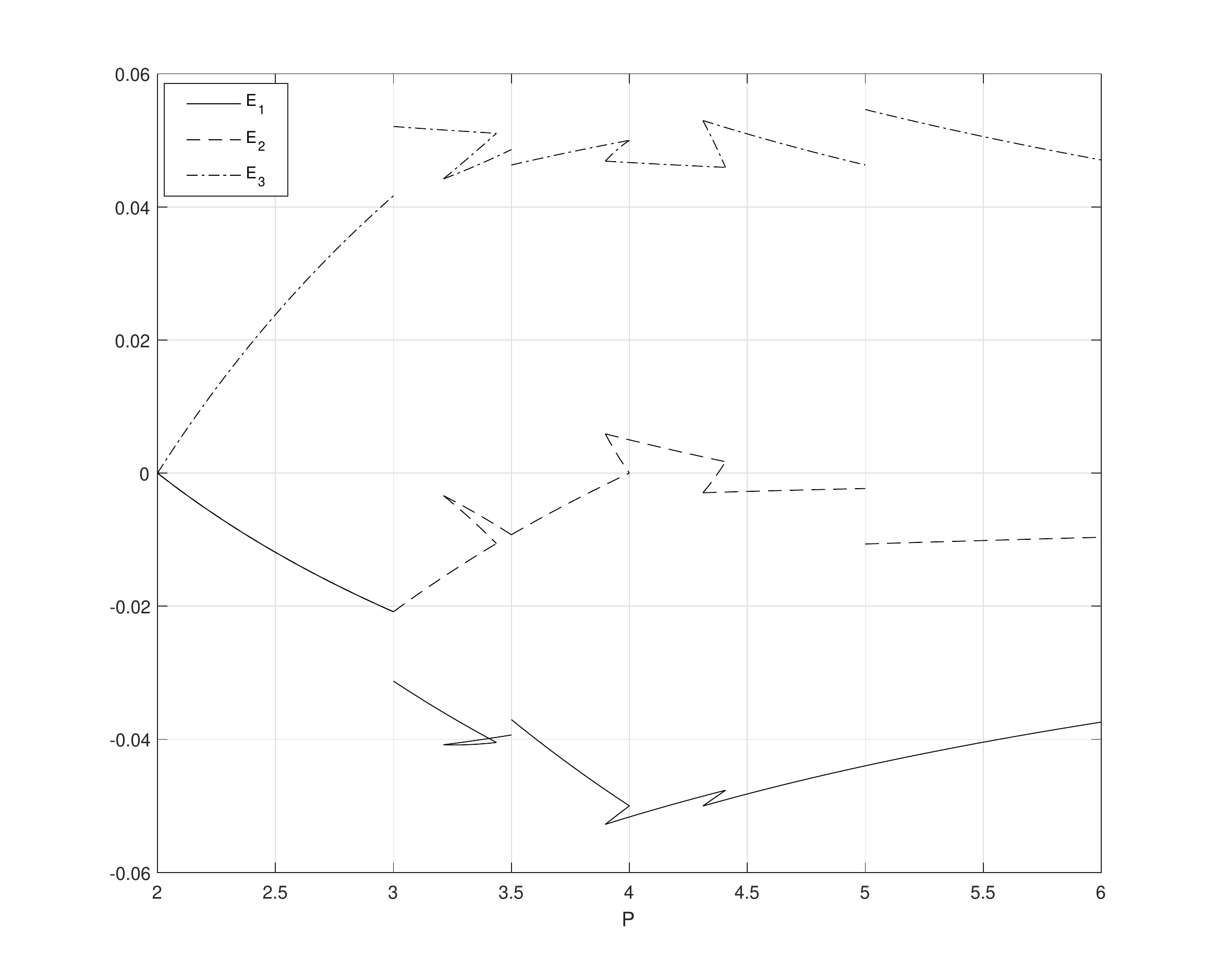}
 \caption{The expected profits at equilibrium.}\label{fig_E}
 \end{center}
 \end{figure}

The expected profits of each player can be calculated using (\ref{eq_E1}) to (\ref{eq_E3}), and are shown in Figure~\ref{fig_E}.
Although the existence of multiple equilibrium solutions makes the structure of these curves quite irregular, Player 3 always has the largest expected profit, because she acts last, whilst for $P>3$, when Player 2 develops a nonzero equilibrium betting frequency, Player 2 has a larger expected profit than Player 1, because she acts after Player 1.

It is also of interest to consider the expected profits for Solutions 2a, 5a and 10a, which are
\begin{equation}
\mbox{\bf Solution 2a ($P=3$):}~24 E_1 = -\frac{1}{2}-\frac{1}{2}a_2,~~24 E_2 = -\frac{1}{2},~~24 E_3 = 1+\frac{1}{2}a_2,~~a_2 \leq \frac{1}{2}.\label{eq_E_sol2a}
\end{equation}
\begin{equation}
\mbox{\bf Solution 5a ($P=\dfrac{7}{2}$):}~24 E_1 = -\frac{10}{9}+\frac{2}{9}a_2,~~24 E_2 = -\frac{2}{9},~~24 E_3 = \frac{4}{3}-\frac{2}{9}a_2~~a_2 \geq \frac{3}{4}.\label{eq_E_sol5a}
\end{equation}
\begin{equation}
\mbox{\bf Solution 10a ($P=5$):}~24 E_1 = -\frac{19}{18},~~24 E_2 = \frac{17}{18}-3b_1,~~24 E_3 = \frac{1}{9}+3b_1,~~\frac{1}{3} \leq b_1 \leq \frac{2}{5}.\label{eq_E_sol10a}
\end{equation}
In Solutions 2a and 5a, Player 2 can transfer profit between Players 1 and 3 at no cost to herself, although it should be noted that for Solution 5a, with $P = \frac{7}{2}$, the equilibrium solution is not unique. In Solution 10a, Player 1 can transfer profit between Players 2 and 3 at no cost to herself. In each case, this can be seen as points of discontinuity in the expected profits plotted in Figure~\ref{fig_E}. This may have implications for repeated play of the game at these pot sizes, as discussed for the one-third street game in \cite{SKP2017} and \cite{Szafron:2013:PFE:2484920.2484962}.

In two player, zero sum games, such as two player Kuhn poker, all equilibrium solutions have the same expected profit. This is not the case for the three player zero sum game that we study here. Moreover, although finding the equilibrium solutions is important, it is the dynamics of the play of the game, which may or may not lead the players towards an equilibrium solution, that determine the outcome of repeated play. 

\section{An ordinary differential equation model of repeated play}\label{sec_odemodels}
In \cite{SKP2017} we studied an ordinary differential equation model of a simplified version of one-third street Kuhn poker. The dynamics are always oscillatory, with only one player bluffcatching, and the identity of that player determined by the initial betting frequencies. We also showed that a more realistic model, which is based on repeated play of the game, has similar dynamics, provided that the players neither adjust their frequencies too rapidly nor have too long a memory. In this section we will study an ordinary differential equation model of the dynamics of repeated play of the simplified full street game. We find that a wider range of dynamics is available than was possible for the system studied in \cite{SKP2017}.

\subsection{The ordinary differential equation model}
We assume that Player $i$ controls the four frequencies $a_i(t)$, $b_i(t)$, $c_i(t)$ and $d_i(t)$ (three frequencies for Player 3, since betting with A is the dominant strategy, which means that $a_3=1$) in order to try to maximise their expected value. For example, we model the behaviour of Player 1's betting frequency, $a_1(t)$, using
\begin{equation}
\frac{da_1}{dt} =24  k_1 a_1 \left(1-a_1\right) \frac{\partial E_1}{\partial a_1} \equiv k_1 a_1 \left(1-a_1\right) \left(-2 b_2 + 2 c_2 - 2 b_3 + 2 d_3 - b_2 c_3\right), \label{eqn_ode1}
\end{equation}
where $k_1$ is a positive constant\footnote{The factor of 24 is purely for computational convenience. It has no effect on the qualitative nature of the dynamics of the system.}. Each frequency $f_i$ satisfies
\begin{equation}
\frac{df_i}{dt} = 24 k_j f_i \left(1-f_i\right) \frac{\partial E_i}{\partial f_i}, \label{eqn_ode2}
\end{equation}
where $k_j$ are positive constants. The term $f_i(1-f_i)$ constrains $f_i$ to lie in $[0,1]$, and also means that, since $df_i/dt=0$ when $f_i = 0$ or $f_i=1$, all of the equilibrium solutions of the game discussed in Section~\ref{sec_fullstreet} are also equilibrium solutions of the system of eleven nonlinear ordinary differential equations defined by (\ref{eqn_ode2}). There also exist equilibrium solutions of the differential equations that are not equilibrium solutions of the game because they do not satisfy one of the constraints on the sign of $\partial E_i/ \partial f_i$ (see Table~\ref{table2}). However, for precisely this reason, these equilibrium solutions are unstable with respect to the dynamics imposed by the differential equations.

From this point onwards, we will fix $k_j = 1$ for $j = 1$, $2$ \ldots $11$ to reduce the complexity of the system. The dynamics of this system can be either periodic, close to periodic or have a long, chaotic transient before being attracted to the boundary of the phase space ($f_i = 0$ or $f_i=1$ for $i = 1$, $2$, $3$), depending on the size, $P$, of the pot. The solution does not asymptote to a steady, stable equilibrium state for any nonequilibrium initial conditions or pot size $P$ for which we solved the equations. We find, through numerical computation of the eigenvalues of each of the equilibrium points given by (\ref{eq_sol1}) to (\ref{eq_sol10}), that Solutions 2, 3, 6, 7 and 8 are unstable. None of Solutions 1, 4, 5, 9 and 10 are asymptotically stable, because each has a centre manifold, but numerical solutions suggest that they are all Lyapunov stable. In other words, solutions that start close to these equilibrium points remain nearby, but do not asymptote to equilibrium as $t \to \infty$. Instead, the solutions are attracted to the centre manifold, where the dynamics are oscillatory. The oscillatory part of the centre manifold has either two dimensions (Solution 1), four dimensions (Solutions 4 and 5), or six dimensions (Solutions 9 and 10). This will become apparent in the nature of the solutions that we discuss below. The increase in the dimensionality of the oscillatory part of the centre manifold reflects the fact that as the pot size, $P$, increases, Player 2, and then for larger $P$, Player 1, has a nonzero bluffing frequency. For $3 < P < P_3 \approx 3.21$ and $4 < P < P_8 \approx 4.31$ all equilibrium points are unstable, and we observe long chaotic transients before the solution asymptotes to the boundary of the phase space, with each betting frequency close to zero or one. Finally, we note that, unlike the ordinary differential equation model studied in \cite{SKP2017}, where the initial betting frequencies determined which of two qualitatively different types of periodic solution emerged as the long time solution, only one type of behaviour is observed for each $P$ that we studied, but with amplitude dependent on the initial betting frequencies.

\subsection{Numerical solution}
\begin{figure}
 \begin{center}
 \includegraphics[width=\textwidth]{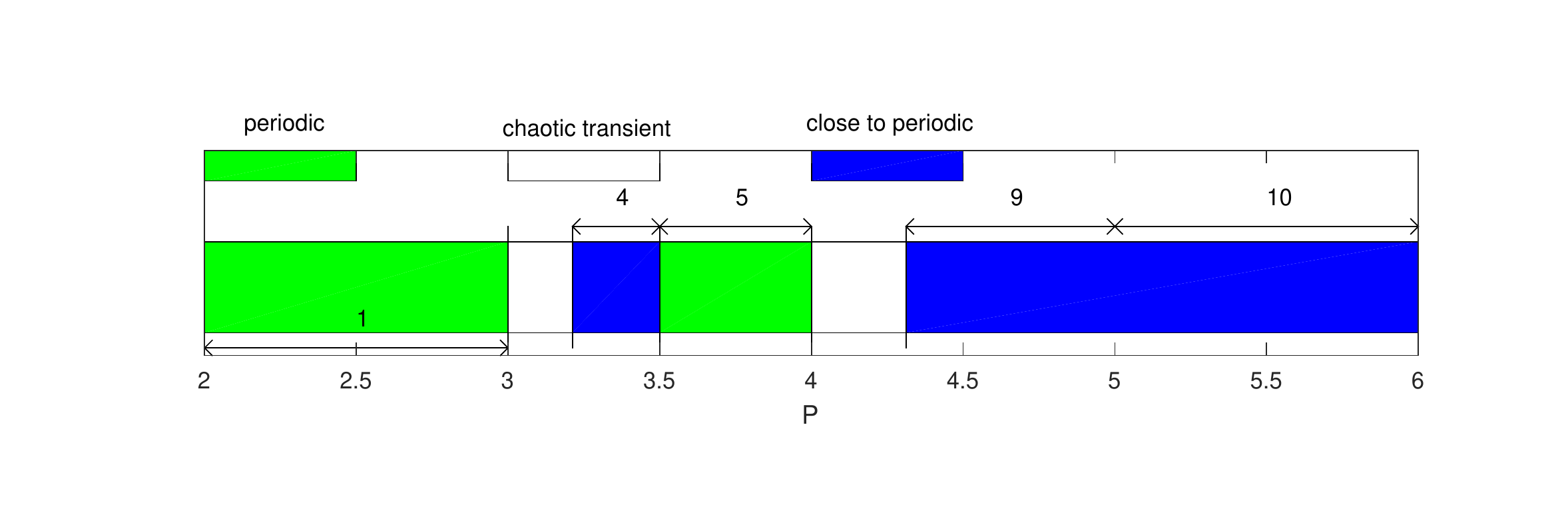}
 \caption{The Lyapunov stable equilibrium solutions, 1, 4, 5, 9 and 10 along with color codes for the behaviour of solutions. Green = periodic; White = chaotic transient; Blue = close to periodic.}\label{fig_ranges2}
 \end{center}
 \end{figure}

In order to solve (\ref{eqn_ode2}) accurately when the betting frequencies are close to zero and one, we solved for the dependent variables $F_i \equiv \log \left\{f_i/\left(1-f_i\right)\right\}$ instead of $f_i$, in terms of which (\ref{eqn_ode2}) take the simpler form
\begin{equation}
\frac{dF_i}{dt} =24 k_j \frac{\partial E_i}{\partial f_i}. \label{eqn_ode3}
\end{equation}
We solved the system of equations (\ref{eqn_ode3}) using the fifth order Runge-Kutta solver {\tt ode45} in MATLAB.  Figure~\ref{fig_ranges2}  shows schematically how the qualitative nature of the dynamics of solutions of (\ref{eqn_ode2}) varies with pot size, $P$. We will give an example of a typical solution from each of these six ranges of values of $P$. We refer to some of the transient behaviour of the solutions as 'chaotic', but note that we have not attempted to study the nature of this in any detail. Our study of this model is only meant to illustrate that the dynamics can be far more complex than those of repeated two player, zero sum games.

\subsubsection{Case 1: $2 < P < 3$}
\begin{figure}
 \begin{center}
 \includegraphics[width=\textwidth]{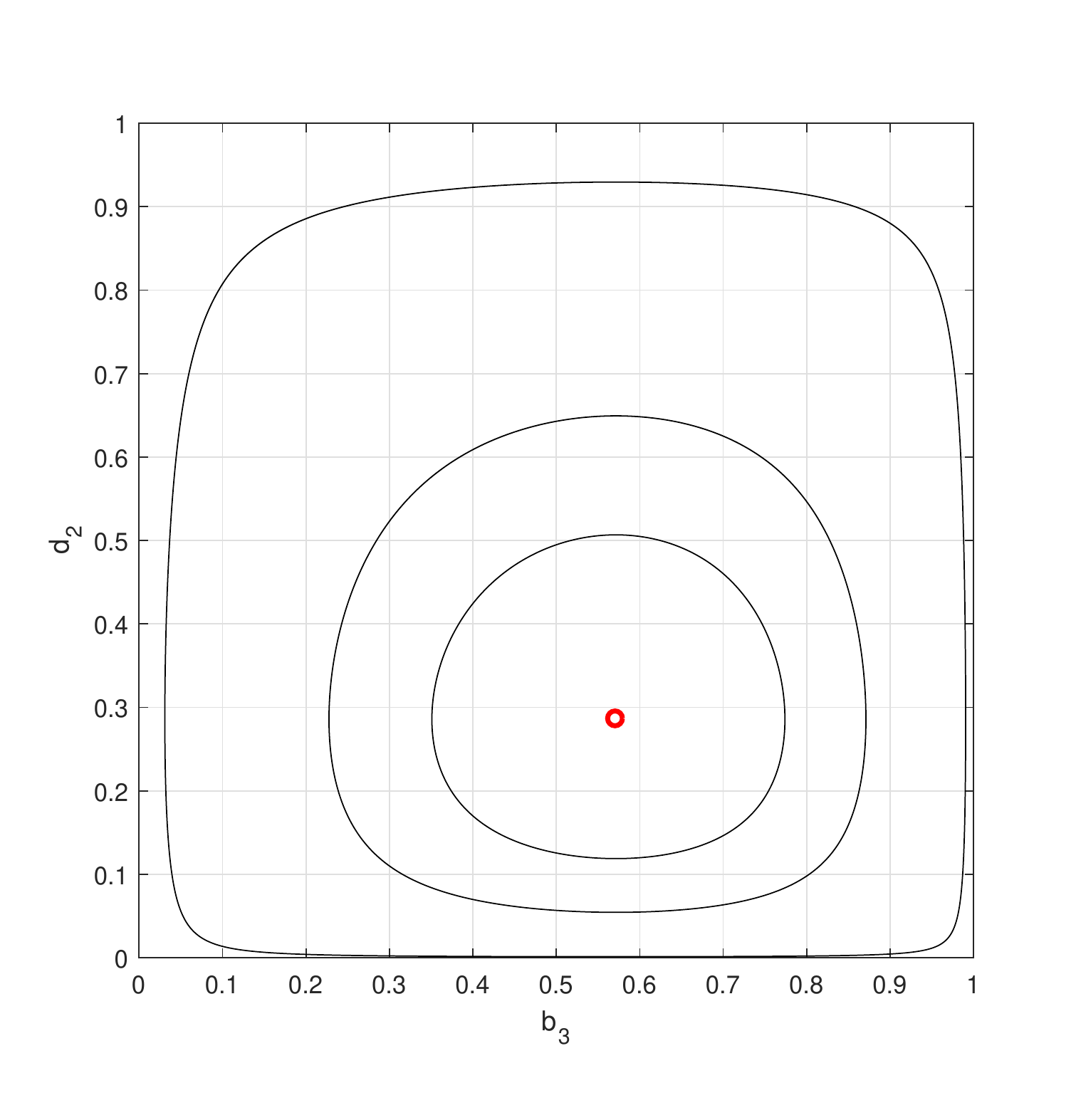}
 \caption{The behaviour of Player 3's bluffing frequency and Player 2's calling frequency for three different sets of initial betting frequencies when $P = 2.5$. Solution 1 is marked with a red circle.}\label{fig_P2p5}
 \end{center}
 \end{figure}

In this case $a_1 \to 0$, $a_2 \to 0$  and $c_1 \to 0$ as $t \to \infty$. The only nontrivial betting frequencies are $b_3$ and $d_2$, which vary periodically with amplitude depending on the initial conditions, as shown in Figure~\ref{fig_P2p5} for $P=2.5$ and three different values of the initial betting frequencies. The centre manifold of Solution 1 attracts the solution. Player 3's bluffing frequency rises and falls, and Player 2, who does the bluff catching, adjusts her calling frequency in response. This behaviour is typical of solutions for $2 < P < 3$ and is the only case for which the dynamics are similar to those of the system studied in \cite{SKP2017}.
\subsection{Case 2: $3<P<P_3 \approx 3.21$}
\begin{figure}
 \begin{center}
 \includegraphics[width=\textwidth]{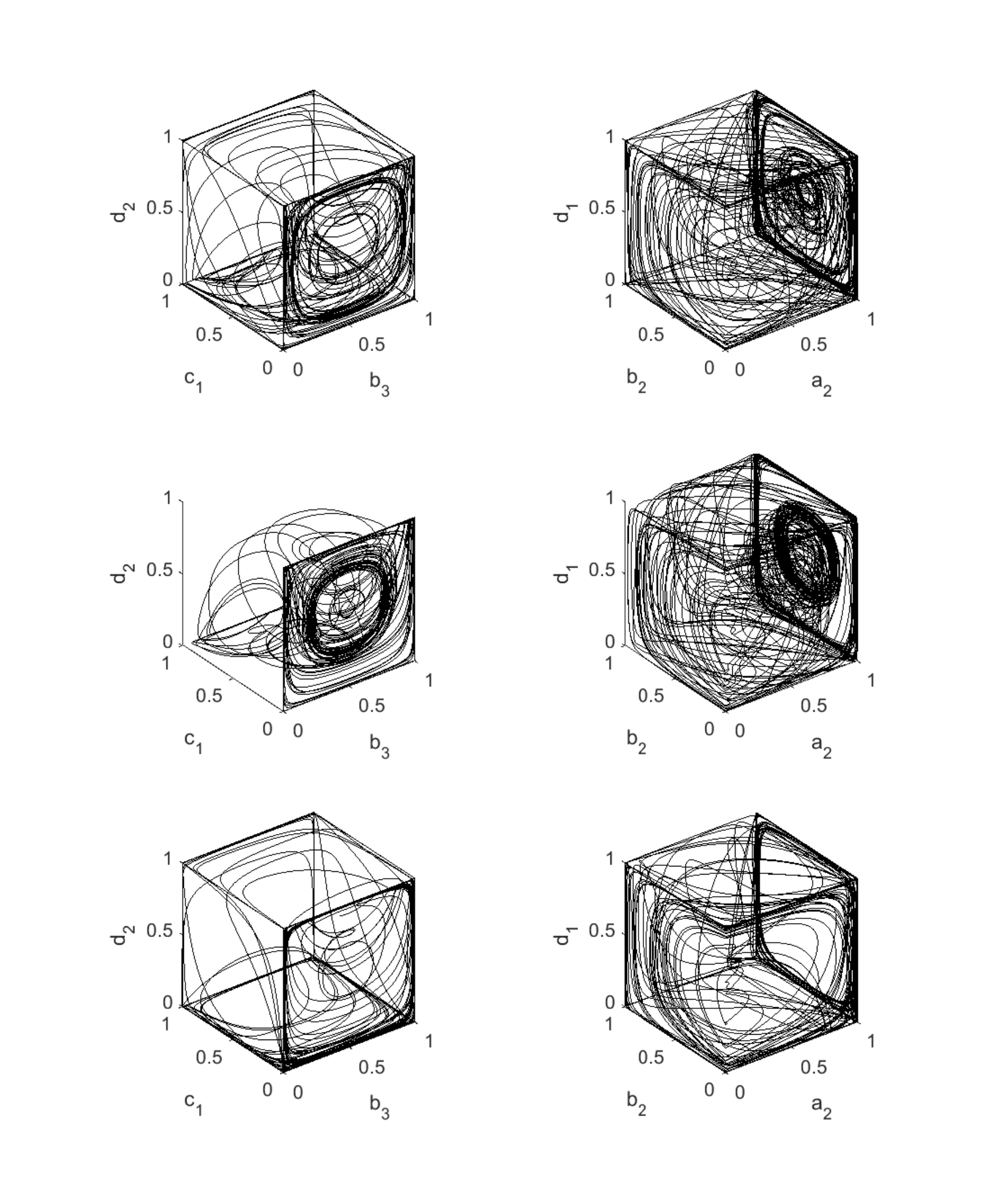}
 \caption{Chaotic transient behaviour when $P=3.1$ for three different initial conditions.}\label{fig_P3p1}
 \end{center}
\end{figure}

In this case, all equilibrium solutions are unstable. Although $a_1 \to 0$, $b_1 \to 0$ and $c_3 \to 0$ as $t \to \infty$, the dynamics of $b_3$, $c_1$, $d_2$, $a_2$, $b_2$ and $d_1$ are complicated, and strongly dependent on the initial betting frequencies, as shown in Figure~\ref{fig_P3p1} for $P=3.1$ and three different values of the initial betting frequencies. Although the initial transient behaviour appears chaotic and typically lasts until $t$ is between $10^3$ and $10^4$, the solution eventually asymptotes to the boundary of the phase space, with each betting frequency close to zero or one. This is likely to be an artefact of the way we have chosen the differential equations (specifically the $f_i(1-f_i)$ terms), and we should simply take it as an indication that the dynamics are nontrivial for this range of pot sizes, for which this behaviour is typical. Similar comments apply in Cases 5 and 6 below.
\subsection{Case 3: $P_3 \approx 3.21 < P < 3.5$}
\begin{figure}
 \begin{center}
 \includegraphics[width=\textwidth]{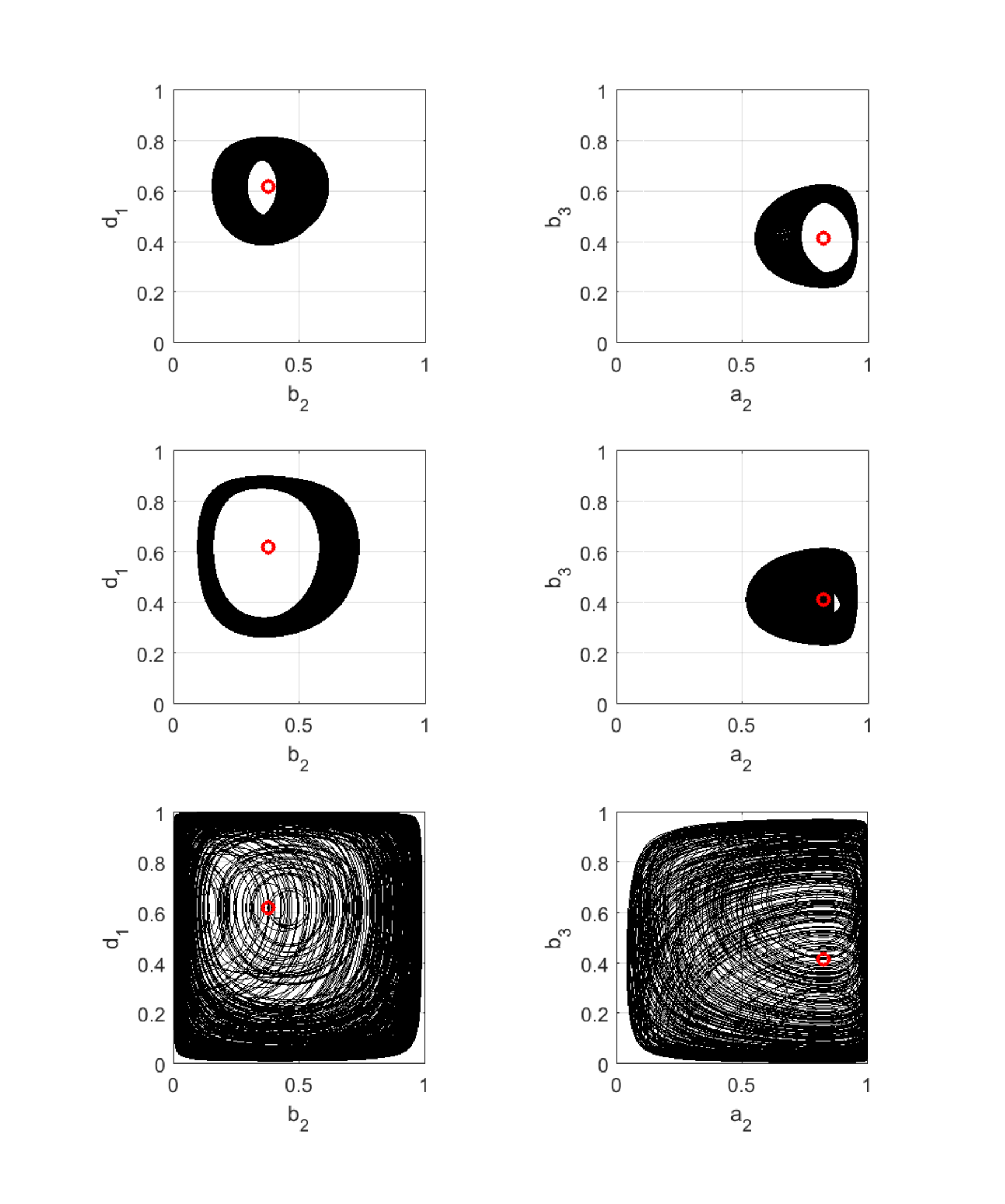}
 \caption{Typical solutions when $P=3.35$ for three different initial conditions. The red circles indicate the location of Solution 4.}\label{fig_P3p35}
 \end{center}
 \end{figure}

\begin{figure}
 \begin{center}
 \includegraphics[width=\textwidth]{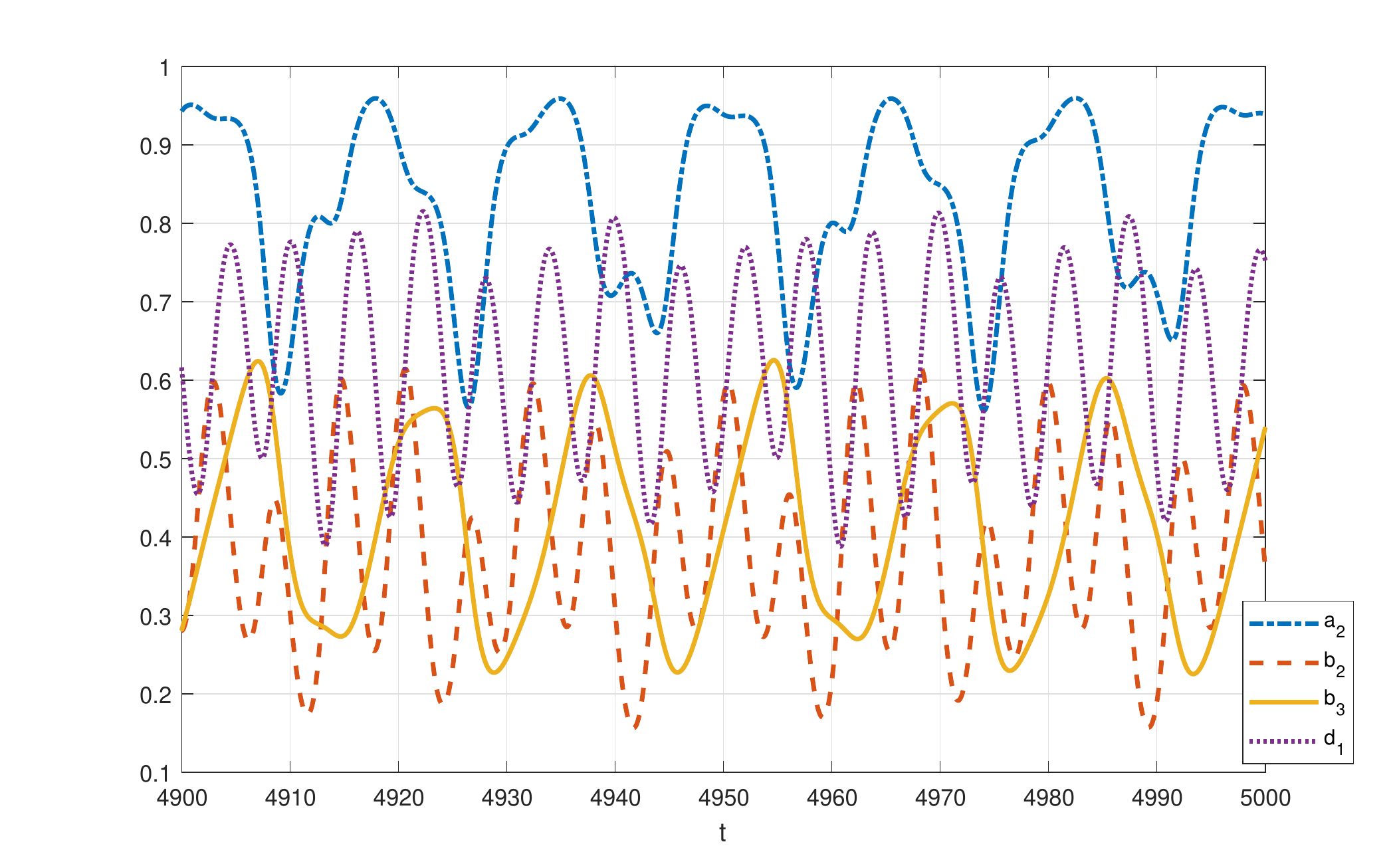}
 \caption{Time series of a typical solution when $P=3.35$.}\label{fig_P3p35_2}
 \end{center}
 \end{figure}

Solution 4 is Lyapunov stable, but not asymptotically stable, and solutions are attracted to its centre manifold, the oscillatory part of which is four-dimensional. The nontrivial betting frequencies are $a_2$, $b_2$, $d_1$ and $b_3$, consistent with the local behaviour of the centre manifold. Some typical solutions are shown in Figures~\ref{fig_P3p35} and~\ref{fig_P3p35_2}. Note that  $a_2$, $b_2$ and $d_1$ are related to adjustments in Player 2's betting and bluffing frequencies and Player 1's associated calling frequency ($c_3 \to 0$). However, $c_1 \to 1$ and $d_2 \to 0$, so the variations in $b_3$ are related indirectly to variations in Player 2's betting frequency $a_2$, which determines how likely she is to have an A with which to call Player 3's bluff. This behaviour is representative of solutions with $P_3 < P < \frac{7}{2}$.
\subsection{Case 4: $3.5<P <4$}
\begin{figure}
 \begin{center}
 \includegraphics[width=\textwidth]{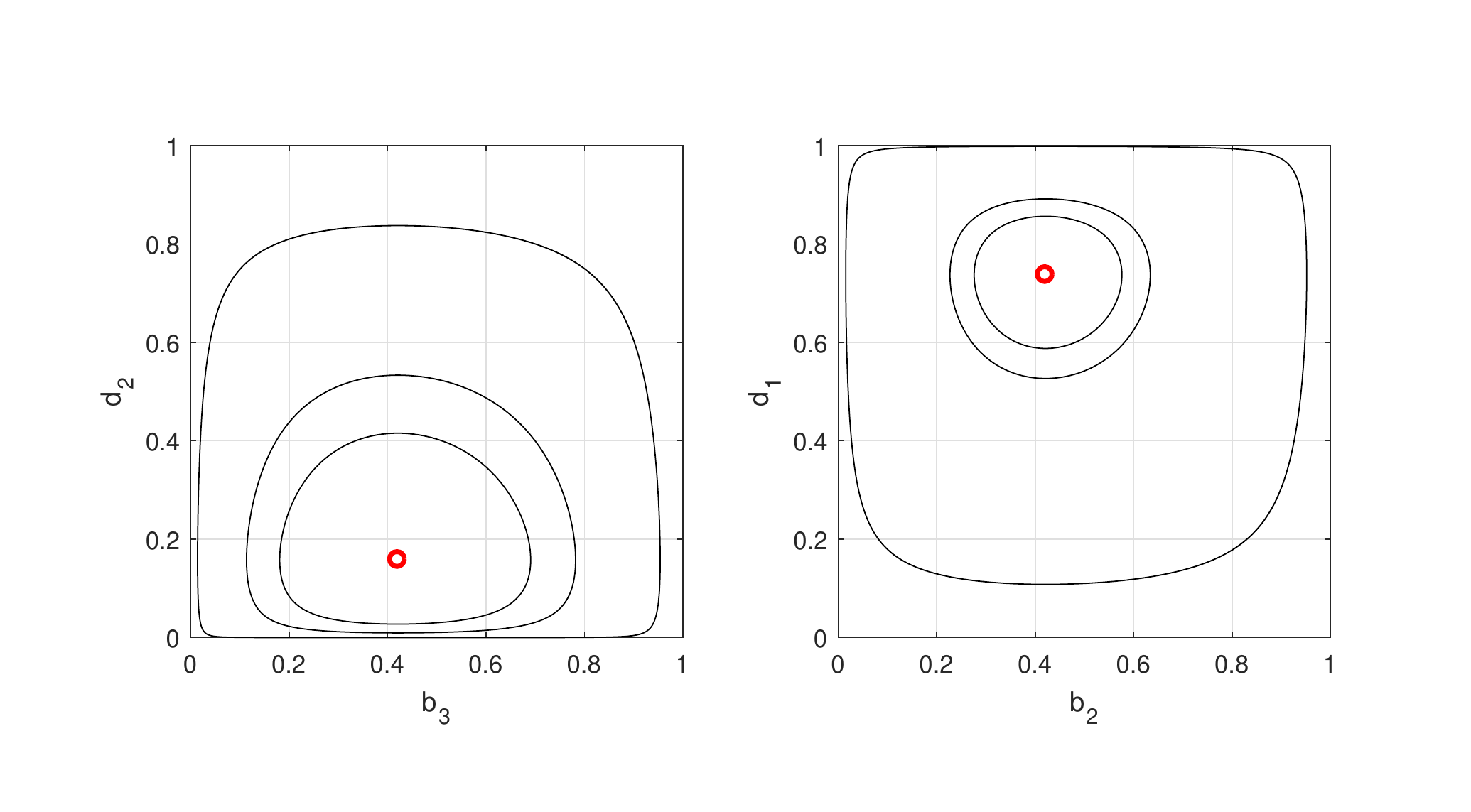}
 \caption{Periodic behaviour when $P=3.75$ for three different initial conditions. The red circles indicate the location of Solution 5.}\label{fig_P3p75}
 \end{center}
 \end{figure}

In this case, as for $P_3<P <3.5$, the solution is attracted to the centre manifold of Solution 5 the oscillatory part of which is four-dimensional, but in this case leads to purely periodic solutions, as shown in Figure~\ref{fig_P3p75}. Note that $a_2 \to 1$, $c_1 \to 1$ and $c_3 \to 0$ as $t \to \infty$. This is typical of solutions with $\frac{7}{2} < P < 4$.
\subsection{Case 5: $4<P<P_8 \approx 4.31$}
\begin{figure}
 \begin{center}
 \includegraphics[width=\textwidth]{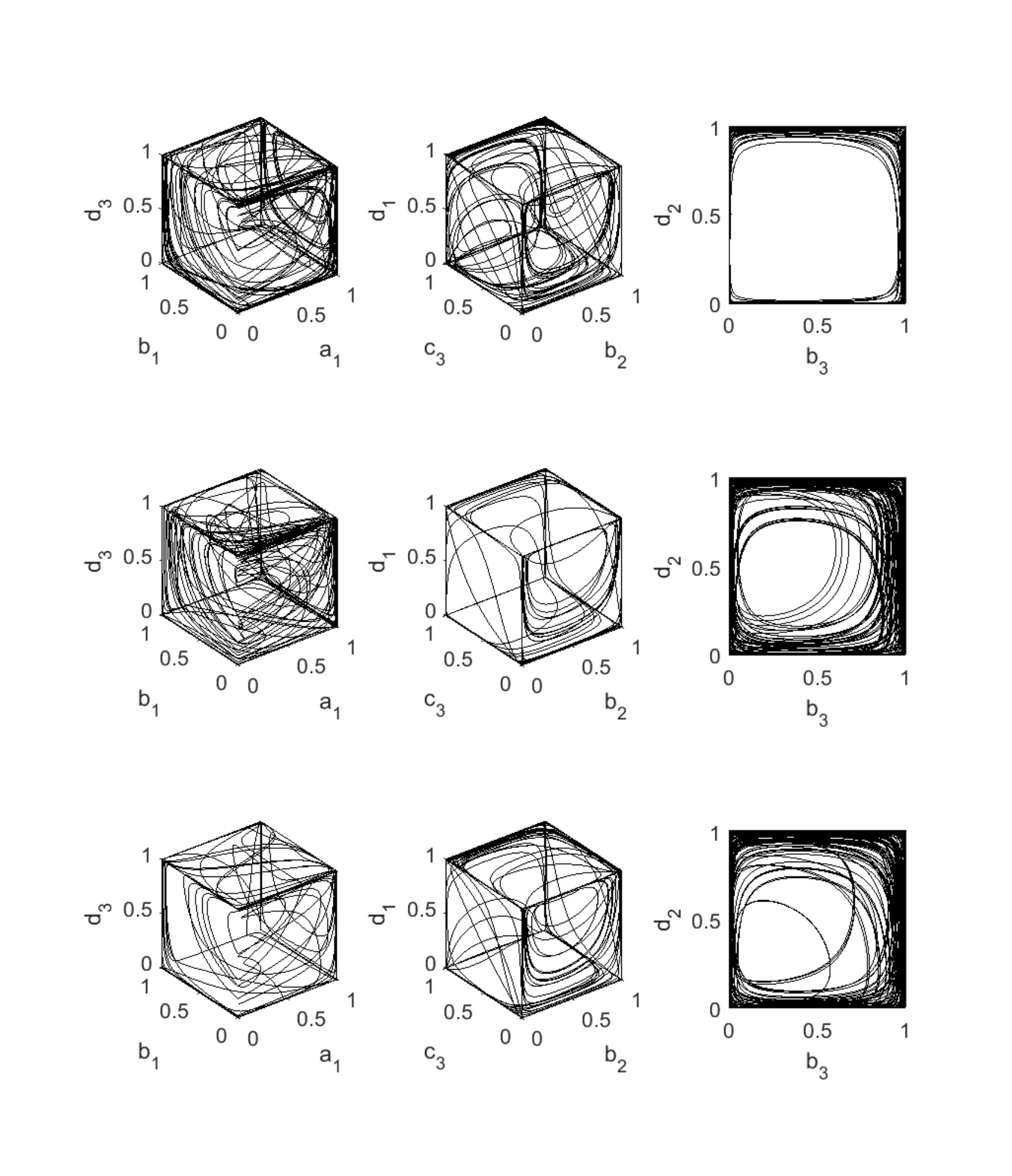}
 \caption{Chaotic transient behaviour when $P=4.15$ for three different initial conditions}\label{fig_P4p15}
 \end{center}
 \end{figure}

In this case, all of the equilibrium solutions are unstable, and the dynamics are similar to those for $3<P<P_3$. Although $c_1 \to 1$, $a_2 \to 1$ and $c_3 \to 0$, there is a complicated initial transient, strongly dependent on the initial betting frequencies, as shown in Figure~\ref{fig_P4p15} for $P=4.15$ and three different values of the initial betting frequencies. However, the solution eventually asymptotes to the boundary of the phase space, with each betting frequency close to zero or one. This behaviour is representative of solutions with $4 < P < P_8$.
\subsection{Case 6: $P>P_8 \approx 4.31$}
\begin{figure}
 \begin{center}
 \includegraphics[width=\textwidth]{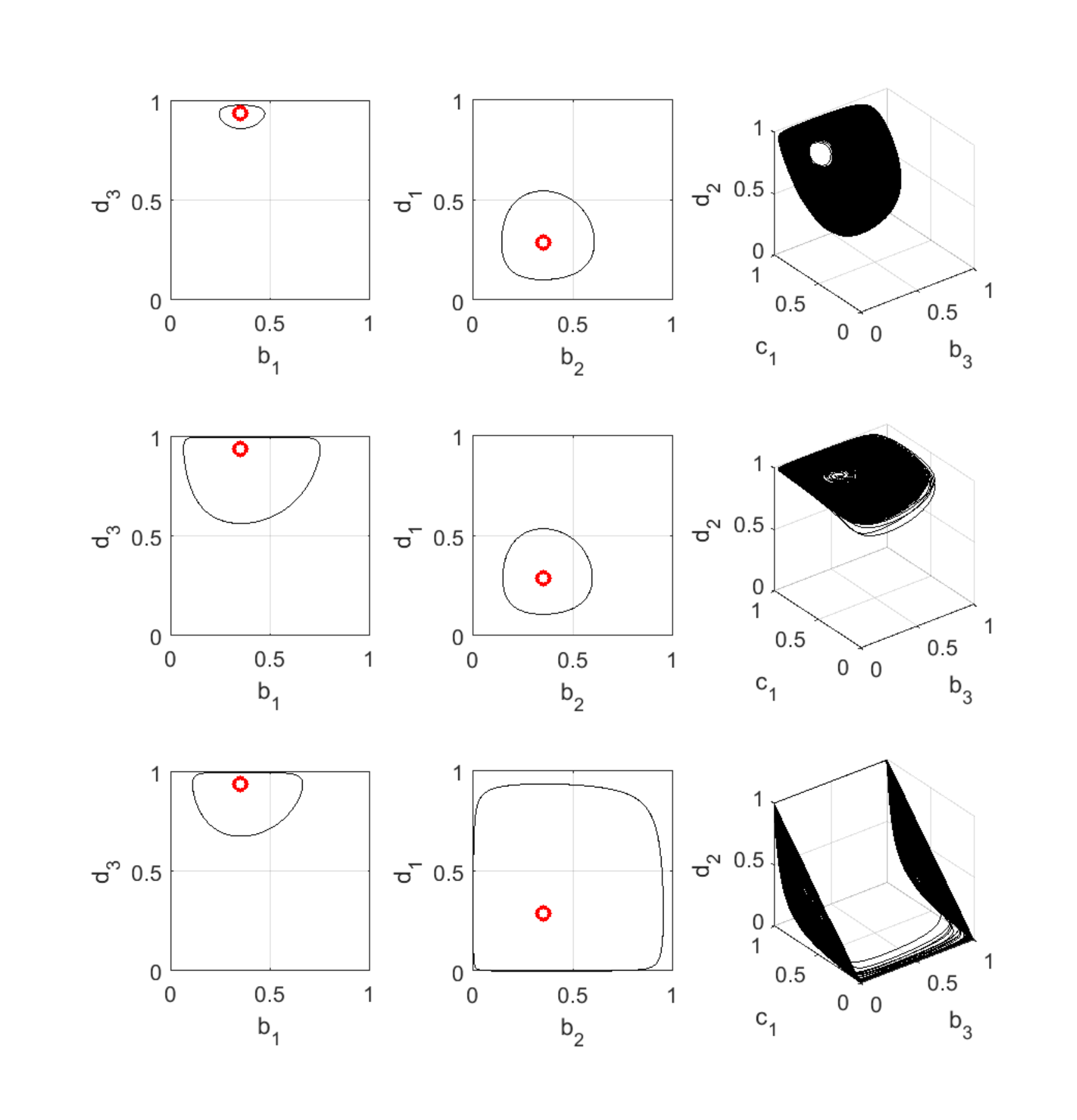}
 \caption{Periodic and close to periodic behaviour when $P=4.65$. The red circles indicate the location of Solution 9.}\label{fig_P4p65}
 \end{center}
 \end{figure}

 In this case, Solution 9 has a seven dimensional centre manifold with six purely  imaginary eigenvalues, and a zero eigenvalue associated with the indeterminacy in $c_1$ and $d_2$ (see (\ref{eq_sol9})). As shown in Figure~\ref{fig_P4p65}, $b_1$ and $d_3$, and $b_2$ and $d_1$ behave periodically. The remaining frequencies, $b_3$, $c_1$ and $d_2$, which are associated with Player 3 bluffing and Players 1 and 2 calling, have a long, close to periodic transient, but are eventually attracted to the boundary of the phase space. For $P>5$, qualitatively similar behaviour occurs around Solution 10.

\begin{figure}
 \begin{center}
 \includegraphics[width=\textwidth]{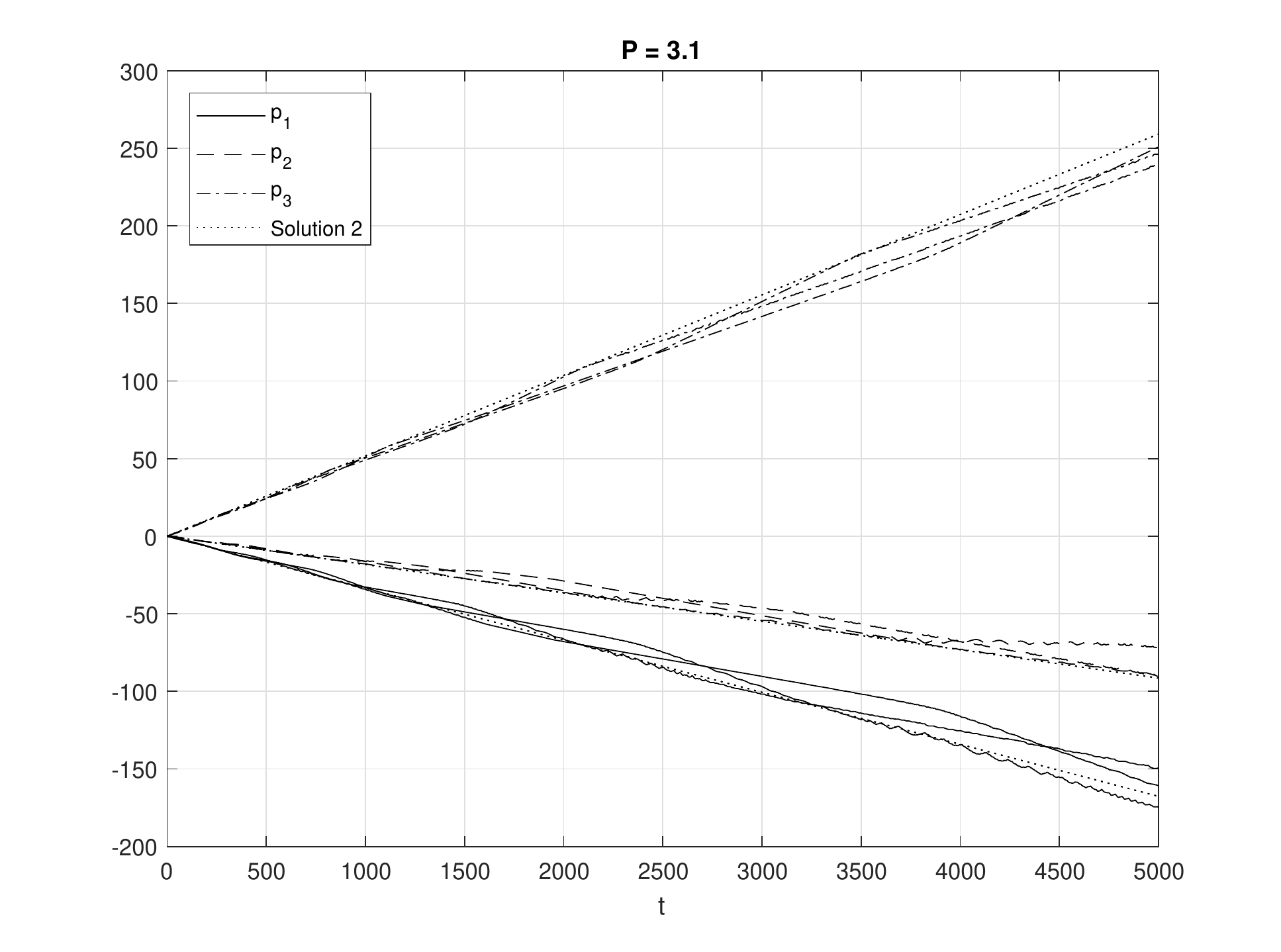}
 \caption{The profits associated with the three solutions for $P = 3.1$ shown in Figure~\protect{\ref{fig_P3p1}}, along with the profit associated with Solution 2.}\label{fig_E1}
 \end{center}
 \end{figure}

In each case, we can also calculate each Player's profit, $p_i$, by solving
\[\frac{dp_i}{dt} = E_i,\]
with the expected profit per game $E_i$ given by (\ref{eq_E1}) to (\ref{eq_E3}). In each of the cases discussed above, each Player's profit remains close to that given by the Lyapunov stable equilibrium solution if one exists, close to Solution 2 for $3 < P <P_3$ and close to Solution 7 for $4 < P < P_8$. Even when the solution has a chaotic transient, as shown in Figure~\ref{fig_P3p1}, the profits $p_i$ shown in Figure~\ref{fig_E1} for these three solutions, closely match the profits associated with Solution 2. 

\section{Conclusions}\label{sec_conc}
In this paper we studied a simplified version of three player Kuhn poker, \cite{Szafron:2013:PFE:2484920.2484962}, in which the card J is only in the deck so that there are more cards than players, making the game nontrivial. Each player must decide whether to bet with A, whether to bluff with Q and whether to catch bluffs with K. This leads to eleven nontrivial betting frequencies.  We found that there are three ranges of values of the pot size, $P$, for which there are three equilibrium solutions. The full structure of the equilibrium solutions was found using the computer algebra package, Mathematica. 

We also studied an ordinary differential equation model of repeated play of the game, in which each player controls their own betting frequencies and tries to increase the rate at which they increase their profit, assuming perfect knowledge of all betting frequencies of every player. Our study of an even simpler game, \cite{SKP2017}, suggests that this gives some insight into the dynamics of repeated play provided that the players adjust their betting frequencies sufficiently slowly, and do not remember the play of the game for too long. We found that the betting frequencies never asymptote to one of the equilibrium solutions, and that the qualitative nature of the solutions depends on the pot size in a manner that is not at all intuitive. Periodic solutions, close to periodic solutions, and solutions with a long chaotic transient are all possible. This illustrates the fact that repeated three player poker games, even games as simple as the one studied here, can have very complex dynamics when the players try to maximise their profit over repeated play of the game, and concepts of balance and equilibrium, much studied by modern professional poker players (referred to as 'GTO' play, \cite{Polk}), may not be as relevant as they think in games with more than two players (but see also \cite{GTORB}).

Finally, we also note that the underlying idea of our ordinary differential equation model, namely that players try to maximise their profits, does not take into account the possibility of cooperation and alliances forming between the players. Although explicit collusion is banned in real world poker games, implicit collusion, in which the players signal to each other within the rules of the game, and no profit is shared outside the game, can occur, and may have a part to play in games such as the one studied in this paper. For an example of how this works in another three player game, see \cite{LSGame2}. 

\bibliography{Kuhn2}{}
\bibliographystyle{hplain}


\end{document}